\documentclass[conference,10pt]{IEEEtran}
\usepackage{amsmath}
\usepackage{resizegather} 
\usepackage{cite}
\usepackage{graphicx}
\usepackage{fixltx2e}
\usepackage{epstopdf} 
\usepackage{array}
\usepackage{amssymb}
\usepackage[english]{babel}
\usepackage{multirow}
\usepackage{booktabs}
\usepackage{pifont}
\usepackage{indentfirst}
\usepackage{algorithm}
\usepackage{algpseudocode}
\usepackage[table]{xcolor}
\usepackage[caption=false,font=footnotesize]{subfig}
\usepackage{placeins}
\usepackage{bm}
\usepackage{array}
\usepackage{float}
\def\checkmark{\tikz\fill[scale=0.4](0,.35) -- (.25,0) -- (1,.7) -- (.25,.15) -- cycle;} 

\addto\captionsenglish{}

\newcommand{\R}{\mathbb{R}}
\usepackage{soul}

\usepackage{cleveref}
\crefname{equation}{}{}
\Crefname{equation}{Equation}{Equations} 
\crefname{table}{Table}{Tables}
\crefname{figure}{Fig.}{Fig.}
\Crefname{figure}{Figure}{Figures}
\crefname{section}{Section}{Sections}

\usepackage{pgfplots}						
\pgfplotsset{compat=newest}
\usepackage{epstopdf}
\usepackage{tikz}
\pgfplotsset{compat=newest}
\pgfplotsset{plot coordinates/math parser=false}
\usepackage{standalone}

\usepackage{fancyhdr}
\usepackage{kantlipsum}

\fancypagestyle{plain}{
  \fancyhf{}
  \fancyhead[L]{Conference on \LaTeX}     
  \fancyfoot[L]{This is a notice}

}
\usepackage{eso-pic}
\begin{document} 
\AddToShipoutPictureBG*{%
  \AtPageUpperLeft{%
    \setlength\unitlength{1in}%
    \hspace*{\dimexpr0.5\paperwidth\relax}
    \makebox(-4.25,-0.75)[c]{\normalsize IEEE PES General Meeting, 2016 (Accepted)}%
}}
\AddToShipoutPictureBG*{%
  \AtPageLowerLeft{%
    \setlength\unitlength{1in}%
    \hspace*{\dimexpr0.5\paperwidth\relax}
    \makebox(-3.20,+1.55)[c]{\footnotesize For better quality of figures, see the published version}%
}}
\title{Quantifying the Effect on the Load Shifting Potential of Buildings due to Ancillary Service Provision}
\author{
	\IEEEauthorblockA{Sarmad Hanif$^{1,2}$
		\IEEEauthorrefmark{1}, Christoph Gruentgens $^{1}$
		\IEEEauthorrefmark{2}, Tobias Massier$^{1}$
		\IEEEauthorrefmark{3}, Thomas Hamacher$^{3}$
		\IEEEauthorrefmark{4}, Thomas Reindl$^{2}$\IEEEauthorrefmark{5}}
	\IEEEauthorblockA{$^{1}$TUM CREATE Limited, Singapore 138602}
	\IEEEauthorblockA{$^{2}$Solar Energy Research Institute of Singapore (SERIS), National University of Singapore (NUS), Singapore 117574}
	\IEEEauthorblockA{$^{3}$Technical University of Munich (TUM), Garching 85748, Germany}
	\IEEEauthorblockA{
		\IEEEauthorrefmark{1}sarmad.hanif@tum-create.edu.sg, 
		\IEEEauthorrefmark{2}christoph.gruentgens@rwth-aachen.de, 
		\IEEEauthorrefmark{3}tobias.massier@tum-create.edu.sg,\\
		\IEEEauthorrefmark{4}thomas.hamacher@tum.de, 
		\IEEEauthorrefmark{5}thomas.reindl@nus.edu.sg}
}
\maketitle

\begin{abstract}
This paper analyzes the pessimistic effect on the inherent load shifting potential (LSP) of buildings due to the participation in the reserve market. A generic model-based optimization approach is deployed, which uses a validated dynamic model along with its experienced external and internal disturbances, to quantify the LSP in the presence of various price signals. The theoretical maximum LSP is obtained using a base energy price without the provision of ancillary services (AS). The deviation from the base case LSP is observed after including the time varying energy and the reserve price from the spot market. Factors affecting the LSP are found to be as: (1) physics of the model, (2) the nature of price signal and (3) the competency of the reserve price with respect to the energy price. Due to its simple formulation, low computation requirements, and modular nature, the method proposed in this paper can easily be deployed by retailers and system operators for the assessment of monetary incentives as well as qualified load type for various demand response (DR) services in liberalized markets.
\end{abstract}



\section{Introduction}
Due to limited controllability of renewable energies, DR services have now been considered all over the world for destressing the grid \cite{DRprog1, DRprog2, EMADR, EMAIL}. Considering the magnitude, buildings are usually considered as one of the major contributors of greenhouse gas emissions and electricity consumption. Moreover, the heating, ventilation and air condition (HVAC) system installed in almost every building is the largest consumer of energy. Hence, the availability of flexibility in the operation of the HVAC system could provide a great potential for DR.

In the past, the developed building/HVAC models \cite{Crawley} served the purpose of estimating annual or monthly energy consumption. Since the need for provision flexibility from the demand-side is now more than ever, building models are required to provide controllable consumption. As a result, a great deal of recent work has been dedicated to developing control-oriented building models \cite{Olde, YMa, Mehdi1}. The idea behind these models is to predict and control the energy consumption of flexible appliances in the building with high fidelity. Most of the control strategies for buildings revolve around the cost optimal consumption \cite{Mehdi2, Mehdi3}, energy efficiency \cite{Mehdi1}, and AS provision to the grid \cite{Vrettosa, Vrettosb, SAR}. The above mentioned literature presents extensive options for the controller design of DR services from buildings. Hence, in order to evaluate these strategies, a generic method must be presented, which quantifies buildings' LSP and the important factors influencing them.

Authors in\cite{Marija} discussed the assessment of types of loads to participate in load management strategies. But the analysis was performed on an unvalidated model. Furthermore, the limits of the LSP and its dependence on AS provision were missing. In \cite{Kochc}, authors provided a method for estimating the storage potential of a cluster of residential thermostatically controlled loads (TCLs). This approach is based upon the hysteresis (ON/OFF) based control methodology of TCLs. Hence, the method is not applicable for a common variable frequency driven fan of an HVAC system. Recently, \cite{KO, RO} provided a method for quantifying demand flexibility. But no comments were obtained regarding the effect to the flexibility with respect to AS provision. Since higher reserve requirements are observed in the presence of renewables, the analysis of building's LSP with the inclusion of AS provision is much needed. 

 Hence, the contribution of this paper is threefold. First, it quantifies the LSP of the modeled building under its physical requirements and monetary incentives. Second, an analytical expression is obtained, explaining the factors related to the AS provision and the deviation in the LSP. Third, the relationship is approximated describing the change in the LSP as a function of energy and reserve price.

Section \ref{sec:1} explains the physical interpretation of the model used for our case study. The method for assessing the LSP and its key findings are given in section \ref{sec:2} and section \ref{sec:3}. Section \ref{sec:4} provides conclusion and possible future research directions related to this topic.
\section{Modeling}\label{sec:1}
The model, which describes the heat transfer of the walls and the air inside the room, consists of five differential equations - four for the temperature of walls and one for the room temperature \cite{Mehdi1}. The nonlinear differential equation for the temperature evolution of the model is given as: 
\begin{equation}\label{eq:1}
	\begin{align}
		\dot{x}_{t} &= f\left(x_t,u_{\mathrm{m},t},u_{\mathrm{d},t},\hat{d_t}\right).
	\end{align}
\end{equation}

Where $x_{t} \in \R^{n}$ represents state vector of size $n = i + j$. $u_{m,t} \in \R^{j}$ and $u_{d,t} \in \R^{j}$ are the input vectors representing the HVAC's mass flow rate and the heat input, respectively. $\hat{d_t}$ is the disturbance vector. The electrical power consumed due to the HVAC's mass flow rate is given as:
\begin{IEEEeqnarray}{LLr}\label{eq:2}
	\IEEEyesnumber\IEEEyessubnumber*
	P_\mathrm{heat} (t)& = \dot{m} (t) c_\mathrm{p} (T_\mathrm{in,heat}(t) - T_\mathrm{room}(t)) \\
	P_\mathrm{fan} (t) &= \frac{\dot{m} (t) \Delta p}{\rho}
\end{IEEEeqnarray}
Where $\dot{m} (t)$, $c_\mathrm{p}$, $T_\mathrm{in,heat} (t)$, $T_\mathrm{room}(t)$, $\rho$ and $\Delta p$ are the HVAC's mass air flow rate, the specific heat capacity of air, the inlet temperature of air, the room temperature, the density of air and the pressure difference across the fan, respectively. The expressions for the consumed power and the differential equation presented above are of nonlinear nature. To be able to control the model with respect to our needs, a linear system is desirable. 

In \cite{Mehdi1}, the author uses Sequential Quadratic Programming (SQP), to obtain a linearized version of \eqref{eq:1}. Zero-order hold is performed to obtain a discretized version of the linear model. The resultant discrete time LTI state space model of the system is given as:
\begin{equation}\label{eq:3}
	\begin{align}
		x_{k+1} &= A x_k + (B_{u} + B_{d}) u_{mk} + E\hat{d}_k
	\end{align}
\end{equation}
Where $x_{k+1} \in \R^{n}$ represents the temperature of all the states $n = i + j$ at step $k+1$ with respect to the inputs $u_\mathrm{m} \in \R^{j}$ and the disturbance $\hat{d}_k \in \R^{n}$ at step k. Matrices $A$, $B_\mathrm{u}$, $B_\mathrm{d}$ and $E$ are of the appropriate sizes.

\begin{figure}[h]
	\centering
	\includegraphics[width=0.5\textwidth]{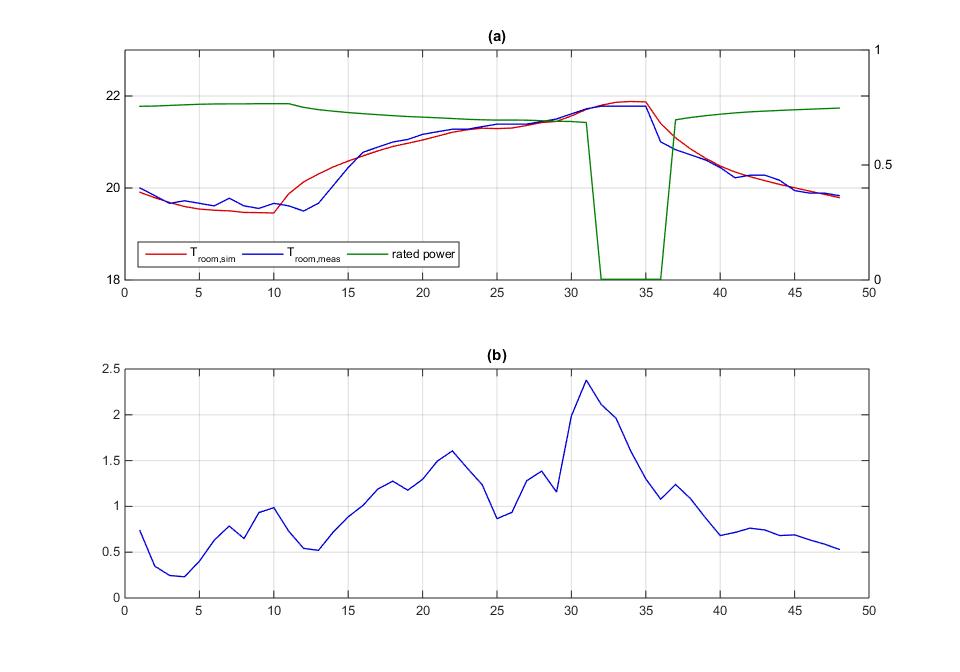}
	\caption{(a) Comparison between measured and simulated room temperature (deg C). The combined rated power of heating and the fan (kW) is shown in green. (b) Disturbance caused by solar radiation on the walls.}
	\label{fig:Temp_comparison_rated_power}
\end{figure}

 The detailed modeling/validation of the system is presented in \cite{Mehdi2}. A simple validation of the model, with our interpretation, is shown in fig. \ref{fig:Temp_comparison_rated_power}. The simulated room temperature from the linear model shows good agreement with the measured temperature. The combined rated power of heating and the fan is shown as the rated power. The drop in the rated power is due to the temperature of the room reaching the maximum temperature set-point (rule-based-controller). The solar disturbance experienced by the room is shown in fig. \ref{fig:Temp_comparison_rated_power}. For further insight into the physical aspects of the model, its parameters, units, and validation, readers are directed to \cite{Mehdi1,Mehdi2,Mehdi3}.

An extension is performed to include the AS provision from the modeled room. Furthermore, to have a sizable provision of flexibility from scheduling of consumption and reserves, another extension is performed to augment the modeled room to a $n_s = 10$ floor building with $n_r = 20$ rooms per floor. The total number of rooms in that building is then $n_{br} = n_s n_r = 200$.
The experienced disturbances for additional rooms are assumed to have 50\% randomly distributed deviation of disturbance from the measured one. With the inclusion of AS provision, state space model of the system is given as:
\begin{equation}\label{eq:4}
	\begin{align}
		\textbf{x}_k  &= \textbf{A} x_{0} + \textbf{B} \textbf{u}_k + \textbf{E} \mathbf{\hat{d}}_{k} + \textbf{B}_{r} \textbf{r}_{k}
	\end{align}
\end{equation}
In \eqref{eq:4}, the state vector for the building $\textbf{x}_k$ and the disturbance $\mathbf{\hat{d}}_{k}$ are of the size $\R^{n n_{br}}$.
 The input vector $\textbf{u}_k \in \R^{n_{br}}$ at step k represents the augmented version of original control variable of the single room $u_\mathrm{m}$.
 $\textbf{B}_{r}$ contains coefficients of the matrix $\textbf{B}$, representing each floor's participation in the AS provision with the quantity $\textbf{r}_{k} \in \R^{n_{br}}$. 
 Matrices $\textbf{A}$, $\textbf{B} = \textbf{B}_u + \textbf{B}_d$, $\textbf{B}_r$ and $\textbf{E}$ are of appropriate sizes.

\section{Quantification of the LSP}\label{sec:2}
\begin{figure}[h]
	\centering
	\includegraphics[width=0.5\textwidth]{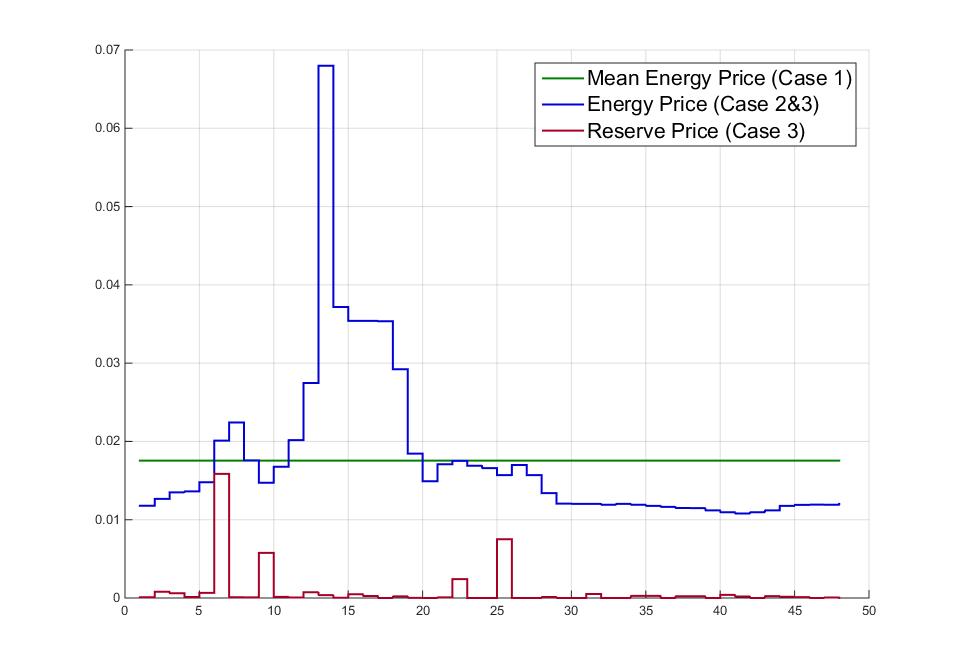}
	\caption{The energy and reserve price in SGD (Singapore Dollar)/kWh taken from the NEMS \cite{NEMS} for the whole day (48 periods - 30 minute steps).}
	\label{EnergyReservePrice_SH}
\end{figure}

To analyze the LSP of the modeled  building in section \ref{sec:1}, sensitivity analysis with respect to change in the electricity price is performed. The reason for opting price-based sensitivity analysis is that the energy price provides an economic equilibrium in the liberalized market consisting of suppliers and consumers. Hence, on the consumption side of the liberalized market of today's smart grid, it is of imperative importance to quantify the price sensitivity of the underlying loads. The LSP quantification procedure is outlined as:
\begin{enumerate}
	\label{cases}
	\item Perform a model-based optimization to obtain a cost optimal pattern of the underlying load subjected to the state feasibility and actuator limits
	\item Perturb the energy price to obtain the response of the load
	\item Fit the desired polynomial function to quantify the response of the load to the change in the price
\end{enumerate}
Seven scenarios are simulated with different perturbations percentiles, as shown in table \ref{table1}. To study the effect of various price signals, each scenario is subjected to three cases, presented in table \ref{table2}. The energy and the reserve price mentioned in table \ref{table1} are shown in fig. \ref{EnergyReservePrice_SH}.
\begin{table}[H]
	\renewcommand{\arraystretch}{1.3}
	\caption{Applied Perturbed Percentiles}
	\label{table1}
	\centering
	\resizebox{\columnwidth}{!}{\begin{tabular}{llllllll}
			\firsthline
			\bfseries Scenario & 1 & 2 & 3 & 4 & 5 & 6 & 7 \\
			\hline
			\bfseries Perturbation  & $\pm10\%$ & $\pm15\%$ & $\pm20\%$ & $\pm25\%$ & $\pm30\%$ & $\pm35\%$ & $\pm40\%$ \\
			\hline
			\lasthline
		\end{tabular}}
	\end{table}
\begin{table}[H]
	\centering
	\caption{Case Studies for each Scenario}
	\label{table2}
	\begin{tabular}{@{}lll@{}}
		\toprule
		Case   & Enegy price          & Reserve price \\ \midrule
		Case 1 & \checkmark (mean)    & --             \\
		Case 2 & \checkmark           & --             \\
		Case 3 & \checkmark           & \checkmark    \\ \bottomrule
	\end{tabular}
\end{table}
Model-based optimization is formulated as: 
\begin{IEEEeqnarray}{LLr}\label{eq:11}
\IEEEyesnumber\IEEEyessubnumber*
\underset{\textbf{u}^{*}, \textbf{r}^{*}}{\text{min}}\ \sum_{k=1}^N \textbf{c}_k^{\bf{T}} \textbf{u}_k - \textbf{k}_k^{\bf{T}} \textbf{r}_k +  \rho \bm{\epsilon}_k \\
\text{subject to} \\
\textbf{x}_{k+1}^{C} = \textbf{A}  \textbf{x}_{k}^{NC} + \textbf{B} \textbf{u}_k + \textbf{E} \mathbf{\hat{d}}_{k} \\
\textbf{x}_{k+1}^{NC} = \textbf{A}  \textbf{x}_{k}^{NC} + \textbf{B} \textbf{u}_k + \textbf{E} \mathbf{\hat{d}}_{k} + \textbf{B}_{r} \textbf{r}_{k}\\
\textbf{x}{^{-}_{k}}  - \bm{\epsilon}_k  \leq  \textbf{x}_{k}^{NC} \leq  \textbf{x}{^{+}_{k}} + \bm{\epsilon}_k\\
\textbf{x}{^{-}_{k}}  - \bm{\epsilon}_k  \leq  \textbf{x}_{k}^{C} \leq  \textbf{x}{^{+}_{k}} + \bm{\epsilon}_k\\
\textbf{u}{^{-}_{k}} \leq \textbf{u}_{k} + \textbf{r}_k \leq \textbf{u}{^{+}_{k}} \label{eq:11con2}\\
\bm{\epsilon}_k, \textbf{u}_{k} - \textbf{r}_k, \textbf{r}_k   \geq \bm{0} \ \ \  \ \ \ \ \ \ \ \ \  \forall k = 0,1 \dots, N 
\end{IEEEeqnarray}
The optimization problem results in the optimal input $\textbf{u}^{*}$ and reserve $\textbf{r}^{*}$ sequence for the whole time duration $N$. In this context, optimality is measured in terms of minimization of the total cost. $\textbf{c}_k $ and $\textbf{k}_k$ convert the consumed energy into the cost of consumption and the revenue from the AS provision, by taking product of the control input and the allocated reserves correspondingly with the energy and the reserve price. The slack variable $\bm{\epsilon}_k \in \R^{nbr}$ in the objective function guarantees the feasibility of the solution by softening the constraints on upper $\textbf{x}^{+}_{k}$ and lower $\textbf{x}^{-}_{k}$ limits of the room. $\rho$ is an arbitrary large term used for penalizing the room temperature. Equation \eqref{eq:11con2} constrains the actuator limits of the HVAC. 

Due to the AS provision, the state is kept feasible for both of the curtailed $\textbf{x}_{k+1}^{C}$ and non-curtailed $\textbf{x}_{k+1}^{NC}$ trajectories. These two trajectories are implemented to replicate the interruptible load program, already in place in the NEMS \cite{EMAIL}. Under this program, load operator can bid a fixed amount of load for each $30$ minute time interval i.e. 1 market period. If the bid is qualified and called upon, the load opeartor must then curtail its load \cite{NEMS}.

The optimization problem presented above is essentially a linear program, and numerous solvers exist which can solve this class of problems very efficiently. For our simulation setup, we have implemented linear programing using YALMIP \cite{YALMIP} and CPLEX \cite{CPLEX}.
\begin{figure}[h]
	\centering
	\includegraphics[width=0.5\textwidth]{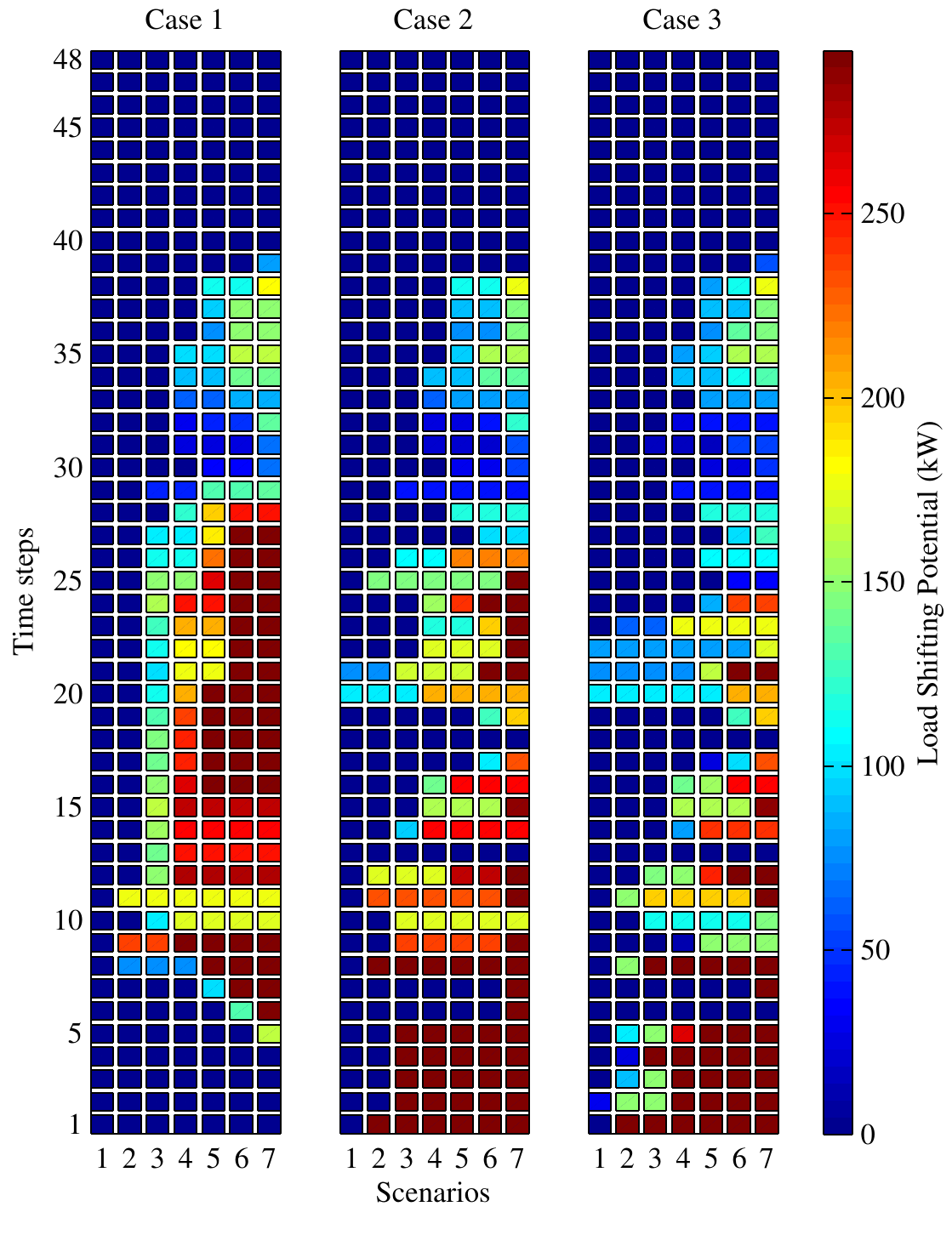}
	\caption{The LSP for all the cases from table \ref{table2}, scenarios from table \ref{table1}, and the price from fig. \ref{EnergyReservePrice_SH}.}
	\label{LoadShift_case123_cropped}
\end{figure}
\section{Simulation Results}\label{sec:3}
\Cref{LoadShift_case123_cropped} shows the LSP of the system. It can be established that due to the load not being shifted uniformly in Case $1$, the thermal inertia of the system is not more than $1$-$2$ time steps. Furthermore, the LSP is shown to be effected by: $(1)$ physical characteristics (heating conditions), $(2)$ comfort requirements (temperature set-point), and $(3)$ energy price. The above mentioned effects are evident from the fact that the LSP is reduced at higher price periods of case $2$ and $3$, and also, beyond the $30^{th}$ time step the system -- due to the lower heating requirements (see fig. \ref{fig:Temp_comparison_rated_power}) -- the LSP is insignificant. Furthermore, for all the cases, it can also be seen that after scenario 4 (percentile $\pm25 \%$), the increase in the LSP ceases. This information is particularly interesting for regulators, responsible for designing incentive schemes for DR services. 

\begin{figure}[h]
	\centering
	\includegraphics[width=0.5\textwidth]{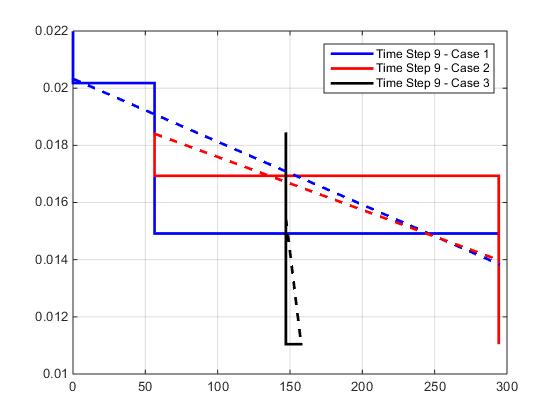}
	\caption{The LSP for all $3$ cases for time step $9$ and scenario $4$ $(\pm25\%)$. y-axis is in SGD/kWh and x-axis in kW. Solid lines represent the optimized consumption for each case and dashed lines show the fitted linear polynomial function.}
	\label{Polyfit_CG}
\end{figure}

\begin{figure}[h]
	\centering
	\includegraphics[width=0.5\textwidth]{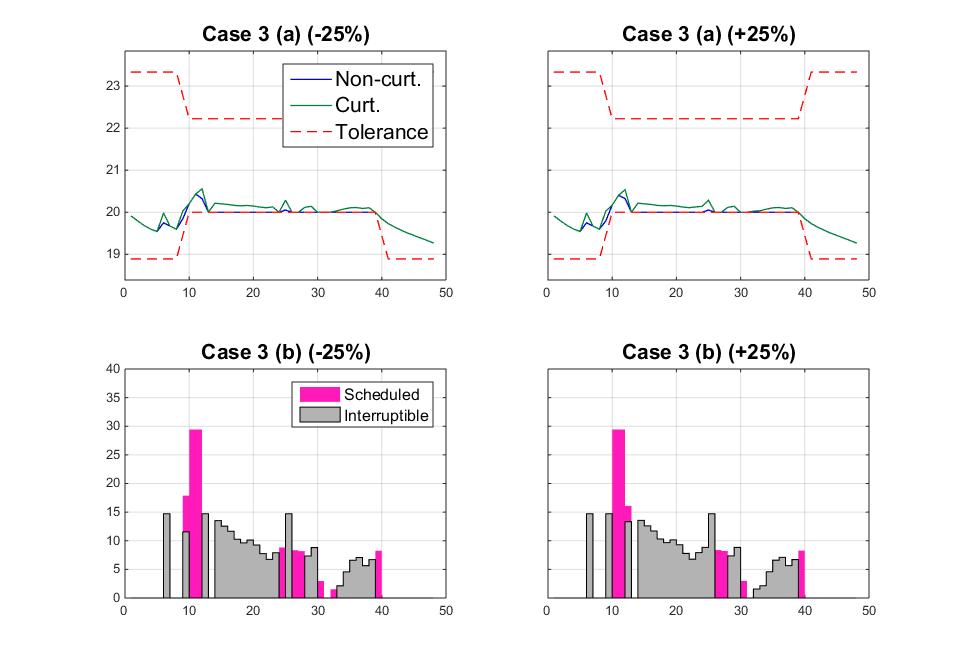}
	\caption{Curtailed and Non-curtailed temperature evolution in deg C(a) and power consumption in kW (b) for the scenario $4$ $(\pm25)$ and the time step $9$.}
	\label{WithReserveCase}
\end{figure}

For scenario $4$, fig. \ref{Polyfit_CG} shows the relationship of the negative effect on the LSP at time step $9$. The polynomial functions in fig. \ref{Polyfit_CG} are approximated as $y(k) = a_1x(k) + a_0$. Where $x(k)$ shows the shiftable demand in kW for the corresponding price $y(k)$ in SGD/kWh. Parameters $a_1$ (slope), and $a_0$ (intercept) are found as explained in the quantification procedure. Physically, $a_1$ and $a_0$ define the elasticity and base value of the shiftable demand, respectively. The first order polynomial functions shown in fig. \ref{Polyfit_CG}, as expected, exhibits the decrease in energy consumption with the increase in the energy price. Compared to the Case $1$ at the time step $9$ for the scenario $4$ $(\pm25)$, the LSP is reduced by ~$17 \%$ and $~96 \%$ for the Case $2$ and $3$, respectively. The demand curves shown in fig. \ref{Polyfit_CG} are of great practical importance, as they quantify the ability to reduce the overall cost of the power system, after being deployed in the well known unit commitment problem of the power system.
\begin{figure}[h]
	\centering
	\includegraphics[width=0.5\textwidth]{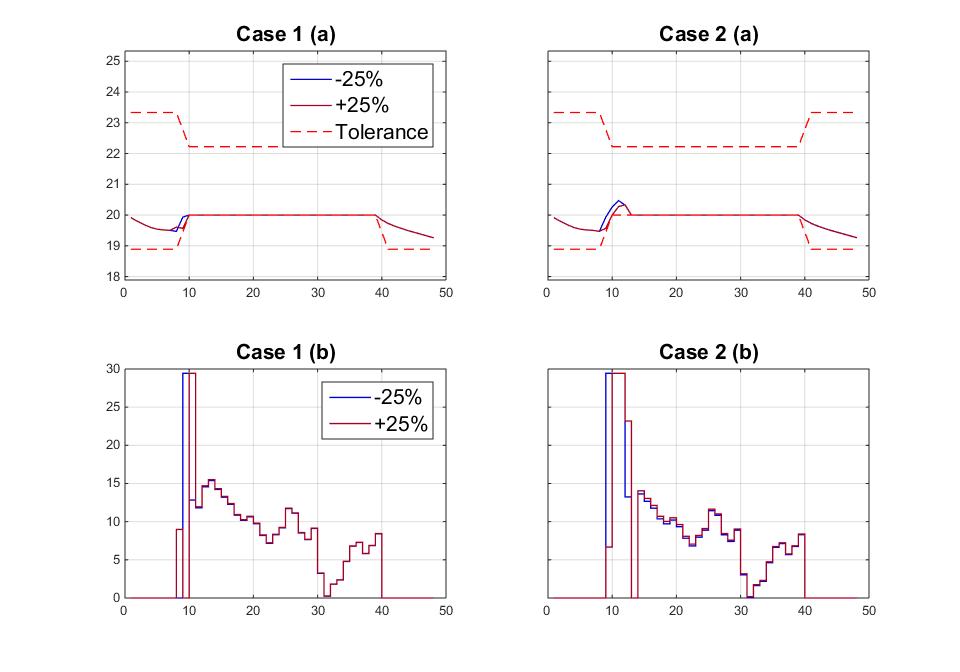}
	\caption{Temperature evolution in deg C(a) and power consumption in kW (b) for scenario $4$ ($\pm25\%$) and time step $9$ .}
	\label{Base_wo_res}
\end{figure}
Providing further insight into the demand curves, fig. \ref{Base_wo_res} and fig. \ref{WithReserveCase} show the temperature evolution of one of the room with its corresponding floor's consumption. In Fig. \ref{Base_wo_res}, compared to Case $1$, the LSP of Case $2$ from $(+25 \%)$ to $(-25 \%)$ is shown to be limited. This is due to the fact that the optimization problem foresees a high price time period (period $12 - 20$ in fig. \ref{EnergyReservePrice_SH}). And, to maintain the comfort requirements of the room, it is not cost optimal anymore to reduce as much load as (compared to Case $1$), even in the event of maximum $(+25 \%)$ perturbation.

The reduction in the LSP is even more pronounced in fig. \ref{WithReserveCase}. The main reason is that now the reserve price at time step $9$ also competes with the energy price of the same period (see fig. \ref{EnergyReservePrice_SH}). Hence, due to the incentives available for AS provision, optimal consumption is not reduced as compared to previous cases of $(+25\%)$ perturbation at time step $9$.
\begin{figure}[h]
	\centering
	\includegraphics[width=0.5\textwidth]{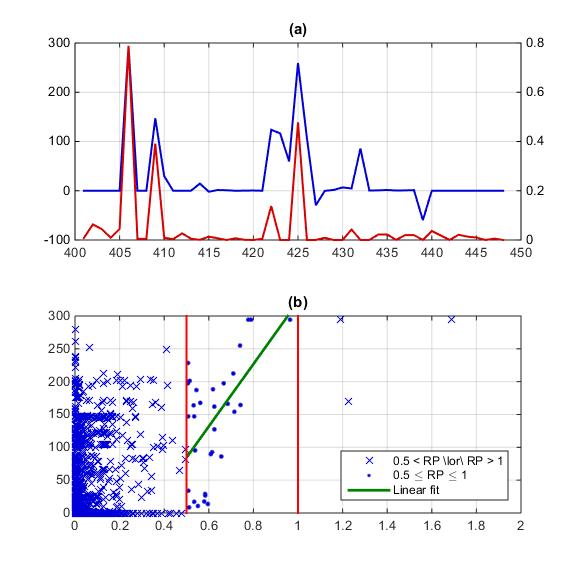}
	\caption{Change in the LSP (kW) with respect to the relative price (RP) (a) for $1$ day and (b) for $6$ months. Region $1$ in (b) is represented as: $0.5\geq$ RP $\leq1$ and region $2$ as: RP $\le0.5$ $\vee$ RP $\ge1$.}
	\label{6month_1day_2in1}
\end{figure}

 To investigate further on the influence of AS provision and the LSP, the procedure outlined in section \ref{sec:2} is simulated for $6$ months of the year $2014$. The subplot (a) in fig. \ref{6month_1day_2in1} is for the time period shown in fig. \ref{LoadShift_case123_cropped}, and the subplot (b) is for the energy and reserve price of $6$ months of the year 2014, taken from the NEMS \cite{NEMS}. The percentage relative price (RP) is the ratio of the reserve price to the energy price. It is evident from the subplot (a) that the increase in the RP deviates the LSP $(\Delta LSP)$ from Case $2$. Subplot (b) has 2 regions. Region $1$ shows an approximate linearly increasing behavior $(\Delta$LSP $=$ $474$RP$ - 153.16)$. Region $2$ does not show any behavior because of two reasons. Firstly, for the values of RP greater than 1, the change in the LSP is saturated due to the physical limits of the modeled system's actuators for the AS provision. Secondly, for values of RP less than $0.5$, the reserve price can not dominate the AS provision term in the optimization problem. This results in the dependence of the change in the LSP for lower values of the RP to be dependent on the factors such as heating and comfort requirements.
\section{Conclusion and Future Work}\label{sec:4}
In this paper, the LSP is quantified and various factors influencing its behavior are explored. Using an optimization based approach, it is shown that the AS provision reduces the natural LSP. Analytical relationships are approximated to explain the demand elasticity. The quantifying procedure explained in this paper could serve to qualify load management or DR schemes, before they are deployed. Future work will include expanding this analysis to the transmission and distribution grid. Furthermore, the effect of technical challenges faced by DR schemes due to network constraints will also be explored.
\section{Acknowledgment}\label{sec:5}
This work was financially supported by the Singapore National Research Foundation under its Campus for Research Excellence And Technological Enterprise (CREATE) programme. This work was also sponsored by National Research Foundation, Prime Minister’s Office, Singapore under its Competitive Research Programme (CRP grant NRF2011NRF-CRP003-030, Power grid stability with an increasing share of intermittent renewables (such as solar PV) in Singapore).
\bibliographystyle{IEEEtran}
\bibliography{ISGTEurope2015_bib}

\begin{thebibliography}{10}
\providecommand{\url}[1]{#1}
\csname url@samestyle\endcsname
\providecommand{\newblock}{\relax}
\providecommand{\bibinfo}[2]{#2}
\providecommand{\BIBentrySTDinterwordspacing}{\spaceskip=0pt\relax}
\providecommand{\BIBentryALTinterwordstretchfactor}{4}
\providecommand{\BIBentryALTinterwordspacing}{\spaceskip=\fontdimen2\font plus
\BIBentryALTinterwordstretchfactor\fontdimen3\font minus
  \fontdimen4\font\relax}
\providecommand{\BIBforeignlanguage}[2]{{%
\expandafter\ifx\csname l@#1\endcsname\relax
\typeout{** WARNING: IEEEtran.bst: No hyphenation pattern has been}%
\typeout{** loaded for the language `#1'. Using the pattern for}%
\typeout{** the default language instead.}%
\else
\language=\csname l@#1\endcsname
\fi
#2}}
\providecommand{\BIBdecl}{\relax}
\BIBdecl

\bibitem{DRprog1}
{J. Torriti, M.G. Hassan and M. Leach.}, ``{Demand response experience in
  Europe: policies, programmes and implementation.}'' in \emph{{Energy}}, vol.
  35(4):.\hskip 1em plus 0.5em minus 0.4em\relax IEEE, 2014, pp. 1575--1583.

\bibitem{DRprog2}
{M. H. Albadi and E. F. El-Saadany}, ``{A summary of demand response in
  electricity markets.}'' in \emph{{Electric Power Systems Research}}, vol.
  78(11):.\hskip 1em plus 0.5em minus 0.4em\relax EPSR, Nov. 2008, pp.
  1989--1996.

\bibitem{EMADR}
\BIBentryALTinterwordspacing
{Energy Market Authority}, ``{Demand Response Programme}.'' [Online].
  Available: \url{http://www.ema.gov.sg/dr/}
\BIBentrySTDinterwordspacing

\bibitem{EMAIL}
S.~Swan, ``{Interruptible load: new partnerships for better energy
  management},'' in \emph{2005 International Power Engineering Conference.},
  vol.~2.\hskip 1em plus 0.5em minus 0.4em\relax IEEE, 2005, pp. 888--892.

\bibitem{Crawley}
{D. B. Crawley. et. al.}, ``{Contrasting the capabilities of building energy
  performance simulation programs.}'' in \emph{{Building and
  Environment.}}\hskip 1em plus 0.5em minus 0.4em\relax Pergamon, 2008.

\bibitem{Olde}
{F. Oldewurtel et. al.}, ``{Energy efficient building and climate control using
  stochastic model predictive control and weather predictions},'' in
  \emph{American Control Conference.}, 2011.

\bibitem{YMa}
{Y. Ma et. al.}, ``{Predictive control for energy efficient buildings with
  thermal storage: Modeling, stimulation, and experiments.}'' in \emph{IEEE
  Control Systems}, vol. 32(1):.\hskip 1em plus 0.5em minus 0.4em\relax IEEE,
  Feb. 2012, pp. 44--64.

\bibitem{Mehdi1}
{M. Maasoumy}, ``{Modeling and optimal conrol algorithm design for HVAC systems
  in energy efficient buildings.}'' in \emph{{EECS Depart. Univ. California,
  Berkeley}}.\hskip 1em plus 0.5em minus 0.4em\relax Master's thesis, Feb.
  2011.

\bibitem{Mehdi2}
{M. Maasoumy et. al.}, ``{Model-based hierarchical optimal control design for
  HVAC systems.}'' in \emph{{Dynamic System Control Conference DSCC.}}\hskip
  1em plus 0.5em minus 0.4em\relax ASME, 2011.

\bibitem{Mehdi3}
{M. Maasoumy and A. S. Vincentelli.}, ``{Total and peak energy consumption
  minimization of building HVAC systems using model predictive control.}'' in
  \emph{{Design and Test of Computers.}}\hskip 1em plus 0.5em minus 0.4em\relax
  IEEE, 2012.

\bibitem{Vrettosa}
{E. Vrettos et. al.}, ``{Predictive control of buildings for demand response
  with dynamic day ahead and real-time prices.}'' in \emph{European Control
  Conference}, June 2013.

\bibitem{Vrettosb}
{E. Vrettos, F. Oldewurtel, F. Zhu and G. Andersson}, ``{Robust Provision of
  Frequency Reserves by Office Building Aggregations.}'' in \emph{19th World
  Congress.}, vol. 19(1):.\hskip 1em plus 0.5em minus 0.4em\relax IFAC., 2014,
  pp. 12\,068--12\,073.

\bibitem{SAR}
{S. Hanif. et. al.}, ``{Model predictive control scheme for investigating
  demand side flexibility in Singapore (In Press).}'' in \emph{{50th University
  Power Engineering Conference.}}\hskip 1em plus 0.5em minus 0.4em\relax IEEE,
  2015.

\bibitem{Marija}
{Marija Ilic' et. al.}, ``{Engineering IT-Enabled Sustainable Electricity
  Services. The tale of Two Low-Cost Green Azores Islands.}'' in
  \emph{{Springer}}.\hskip 1em plus 0.5em minus 0.4em\relax Springer, 2013.

\bibitem{Kochc}
{S. Koch et. al.}, ``{Active Coordination of Thermal Household Appliances for
  Load Management Purposes.}'' in \emph{{Symposium on Power Plants and Power
  Systems Control.}}\hskip 1em plus 0.5em minus 0.4em\relax IFAC, {5-8 July}
  2009.

\bibitem{KO}
{K. Kouzelis et. al.}, ``{Probabilistic quantification of potentially flexible
  residential demand.}'' in \emph{{PES}}.\hskip 1em plus 0.5em minus
  0.4em\relax IEEE, 2014.

\bibitem{RO}
{Roel De Coninck, L. Helsen}, ``{Bottom-up quantification of the flexibility
  potential of buildings.}'' in \emph{{13th International Conference of the
  International Building Performance Simulation Association.}}\hskip 1em plus
  0.5em minus 0.4em\relax Building Simulation, 2013.

\bibitem{NEMS}
{F. Lu and Henry Gan}, ``{National Electricity Market of Singapore},'' 2005.

\bibitem{YALMIP}
{J. Loefberg.}, ``{YALMIP : A toolbox for modeling and optimization in
  MATLAB.}'' in \emph{{International Symposium on Computer Aided Control System
  Design}}.\hskip 1em plus 0.5em minus 0.4em\relax IEEE, 2004.

\bibitem{CPLEX}
\BIBentryALTinterwordspacing
{IBM}, ``{IBM ILOG CPLEX Optimization Studio}.'' [Online]. Available:
  \url{{http://www-03.ibm.com/software/products/en/ibmilogcpleoptistud}}
\BIBentrySTDinterwordspacing

\end{thebibliography}
\end{document}